\documentclass[preprint,11pt]{elsarticle}
\usepackage{amssymb, pstricks,pst-plot}
\journal{Comptes Rendus Mathematique}

\input xy
\xyoption{all}

\parskip=\smallskipamount

\newtheorem{theorem}{Theorem}[section]

\newtheorem{corollary}[theorem]{Corollary}
\newtheorem{proposition}[theorem]{Proposition}

\newtheorem{question}[theorem]{Question}

\newcommand{\C}{\mathbb{C}}

\newcommand{\R}{\mathbb{R}}

\newcommand{\cE}{\mathcal{E}}

\newcommand{\cO}{\mathcal{O}}

\newcommand\wt{\widetilde}

\begin{document}

\begin{frontmatter}
\title{Oka Maps} 
\author{Franc Forstneri\v c}
\address{Faculty of Mathematics and Physics, University of Ljubljana, 
and Institute of Mathematics, Physics and Mechanics, Jadranska 19, 
1000 Ljubljana, Slovenia}

\begin{abstract}
We prove that for a holomorphic submersion of reduced complex spaces, 
the basic Oka property implies the parametric Oka property. 
It follows that a stratified subelliptic submersion, or a 
stratified fiber bundle whose fibers are Oka manifolds, 
enjoys the parametric Oka property. 
\smallskip

\noindent 
{\bf R\'esum\'e}
\smallskip

\bf L\'es Applications d'Oka. \rm
Nous prouvons que, pour une submersion holomorphe des espaces complexes
r\`eduits, la propri\'et\'e d'Oka simple implique la propri\'et\'e d'Oka 
param\'etrique. En particulier, toute submersion sous-elliptique 
stratifi\'e poss\`ede la propri\'et\'e d'Oka param\'etrique. 
\end{abstract}

\end{frontmatter}

%
%
%
%   OKA PROPERTIES OF HOLOMORPHIC MAPS
%
%
%
\section{Oka properties of holomorphic maps}
\label{Okaproperties}
Let $E$ and $B$ be reduced complex spaces.
A holomorphic map $\pi\colon E\to B$ is said to enjoy the 
{\em Basic Oka Property} (BOP) if, given a 
holomorphic map $f\colon X\to B$ from a reduced Stein space $X$ 
and a  continuous map $F_0\colon X\to E$ 
satisfying $\pi\circ F_0=f$ (a {\em lifting} of $f$)
such that $F_0$ is holomorphic on a closed complex subvariety
$X'$ of $X$ and in a neighborhood of a compact 
$\cO(X)$-convex subset $K$ of $X$, there is a 
homotopy of liftings $F_t\colon X\to E$ $(t\in[0,1])$ of $f$
to a holomorphic lifting $F_1$ such that 
for every $t\in[0,1]$, $F_t$ is holomorphic in a neighborhood 
of $K$ (independent of $t$), 
$\sup_{x\in K} {\rm dist}(F_t(x),F_0(x))<\epsilon$, and 
$F_t|_{X'} =F_0|_{X'}$ (the homotopy is fixed on $X'$).

% Here ${\rm dist}$ denotes a distance function on $E$
% which induces the underlying complex space topology.

By definition, a complex manifold $Y$ enjoys BOP 
if and only if the trivial map $Y\to point$ does.
This is equivalent to several other properties, from the simplest 
{\em Convex Approximation Property} (CAP) 
to the {\em Parametric Oka Property} (POP) concerning compact 
families of maps from reduced Stein spaces to $Y$ \cite{FF:CR}. 
A complex manifold enjoying these equivalent properties is called 
an {\em Oka manifold} \cite{FF:CR,Larusson4}; 
these are precisely the {\em fibrant complex manifolds} 
in L\'arusson's model category \cite{Larusson2}.
Here we prove that $\hbox{\rm BOP}\Longrightarrow \hbox{\rm POP}$ 
also holds for holomorphic submersions.
(The submersion condition corresponds 
to requiring smoothness as part of the definition of a variety 
being Oka.  The singular case is rather problematic.)

\begin{theorem}
\label{BOP-POP}
For every holomorphic submersion $\pi\colon E\to B$ of reduced complex spaces, 
the basic Oka property implies the parametric Oka property.
\end{theorem}

Recall \cite{Larusson2} that a holomorphic map $\pi\colon E\to B$ 
enjoys the {\em Parametric Oka Property} (POP) if
for any triple $(X,X',K)$ as above and for any pair $P_0\subset P$ 
of compact subsets in a Euclidean space $\R^m$ the following holds. 
Given a continuous map $f\colon P\times X\to B$ that is {\em $X$-holomorphic}
(that is, $f(p,\cdotp)\colon X\to B$ is holomorphic for every $p\in P$)
and a continuous map $F_0 \colon P\times X\to E$ such that (a) $\pi\circ F_0=f$,
(b) $F_0(p,\cdotp)$ is holomorphic on $X$ for all $p\in P_0$
and is holomorphic on $K\cup X'$ for all $p\in P$, there 
exists for every $\epsilon >0$ a homotopy of continuous liftings 
$F_t\colon P\times X\to E$ of $f$ to an $X$-holomorphic lifting
$F_1$ such that the following hold for all $t\in[0,1]$:
\begin{itemize}
\item[\rm (i)]   $F_t=F_0$ on $(P_0\times X) \cup (P\times X')$, and
\item[\rm (ii)]  $F_t$ is $X$-holomorphic on $K$ and 
$
	\sup_{p\in P,\, x\in K} {\rm dist}(F_t(p,x),F_0(p,x))<\epsilon. 
$
\end{itemize}

A stratified subelliptic holomorphic submersion, 
or a stratified fiber bundle with Oka fibers, enjoys BOP 
\cite{FF:Kohn,FF:Rothschild}. Hence Theorem \ref{BOP-POP} 
implies:

\begin{corollary}
\label{SES-POP}
{\rm (i)} Every stratified subelliptic submersion enjoys {\em POP}. 

\noindent {\rm (ii)} 
Every stratified holomorphic fiber bundle with Oka fibers enjoys {\em POP}. 
\end{corollary}

If $\pi\colon E \to B$ enjoys the Oka property then by 
considering liftings of constant maps $X\to b\in B$ 
we see that every fiber $E_b=\pi^{-1}(b)$ is an Oka manifold.
For stratified fiber bundles the converse holds 
by Corollary \ref{SES-POP}.

\begin{question}
Does every holomorphic submersion
with Oka fibers enjoys the Oka property?
\end{question}

A holomorphic map is said to be an {\em Oka map} 
if it is a topological (Serre) fibration and it enjoys POP. 
Such maps are {\em intermediate fibrations} 
in L\'arusson's model category \cite{Larusson2,Larusson3}. 
Corollary \ref{SES-POP} implies:

\begin{corollary}
\label{cor:Okamaps}
(i) Every holomorphic fiber bundle projection with Oka fiber is an Oka map.

(ii) A stratified subelliptic submersion, or a stratified 
holomorphic fiber bundle with Oka fibers, is an Oka map 
if and only if it is a Serre fibration.
\end{corollary}

Corollary \ref{SES-POP} (i) and the proof by 
Ivarsson and Kutzschebauch \cite{IK} gives the following 
solution of the parametric Gromov-Vaser\-stein problem
\cite{Gromov:OP,Vaser}.

\begin{theorem}
\label{Ivar-Kut} 
Assume that $X$ is a finite dimensional reduced Stein space, 
$P$ is a compact subset of $\R^m$, and  
$f\colon P\times X\to {\rm SL}_n(\mathbb{C})$ is a null-homotopic 
$X$-holomorphic mapping. Then there exist a natural number $N$ and 
$X$-holomorphic mappings 
$G_1,\dots, G_{N}\colon P\times X \to \mathbb{C}^{n(n-1)/2}$ such that
\[
		f(p,x) = \left(\begin{array}{cc} 1  &  0 \cr G_1(p,x) & 1 \cr \end{array} \right)  
		\left(\begin{array}{cc} 1 & G_2(p,x) \cr 0 & 1 \cr \end{array} \right)  \ldots 
		\left(\begin{array}{cc} 1 & G_N(p,x)\cr 0 & 1 \cr \end{array} \right)
\]
is a product of upper and lower diagonal unipotent matrices.
\end{theorem}

%
%
%
%  REDUCTIONS
%
%
\section{Reduction of Theorem \ref{BOP-POP} to an approximation property}
\label{Reductions}
Assume that $\pi\colon E\to B$ enjoys BOP and that
$(X,X',K,P,P_0,f,F_0)$ are as in the definition of POP,
with $P_0\subset P\subset \R^m\subset \C^m$. Set
\begin{equation}
\label{eq:psi}
	Z=\C^m \times X\times E,\quad Z_0 = \C^m \times X\times B,
	\quad \psi = ({\rm id}_{\C^m \times X}) \times \pi \colon Z\to Z_0.
\end{equation}
Observe that $\psi$ enjoys BOP (resp.\ POP) if and only if $\pi$ does.
To the map $f\colon P\times X \to B$ we associate the 
$X$-holomorphic section 
\begin{equation}
\label{eq:g}
	g \colon P\times X\to Z_0,\quad 
  g(p,x)=(p,x,f(p,x)) \quad (p\in P,\ x\in X),
\end{equation}
and to the $\pi$-lifting $F_0\colon P\times X\to E$ of $f$ 
we associate the section 
\begin{equation}
\label{eq:G}
	G_0 \colon P\times X\to Z,\quad 
	G_0(p,x)=\left( p,x,F_0(p,x)\right) \quad (p\in P,\ x\in X).
\end{equation}
Then $\psi \circ G_0 = g$, $G_0$ is $X$-holomorphic over 
$K\cup X'$, and $G_0|_{P_0\times X}$ is $X$-holomorphic.
We must find a homotopy $G_t\colon P\times X\to Z$ 
$(t\in[0,1])$ such that $\psi \circ G_t=g$ for all $t\in [0,1]$,
$G_1$ is $X$-holomorphic, and for all $t\in [0,1]$ 
the map $G_t$ has the same properties as $G_0$,  
$G_t$ is uniformly close to $G_0$ on $K\times P$, and
$G_t=G_0$ on $(P_0\times X)\cup (P\times X')$. Set 
\[
		Q= [0,1] \times P, \quad Q_0=(\{0\} \times P)\cup ([0,1]\times P_0).
\]
The following result is the key to the proof of Theorem \ref{BOP-POP}. 

%The set $K$ in Proposition \ref{BOP-PHAP} is not necessarily the same as above.

%
%
%  BOP IMPLIES PHAP
%
%
\begin{proposition}
\label{BOP-PHAP}
If the submersion $\psi\colon Z\to Z_0$ {\rm (\ref{eq:psi})} 
enjoys the basic Oka property, then it also enjoys the following 

\noindent {\em Parametric Homotopy Approximation Property (PHAP)}:
Let $K\subset L$ be compact $\cO(X)$-convex subsets
and let $U\supset K$, $V\supset L$ be open neighborhoods
in $X$. 
% Let ${\rm dist}$ be a distance on $Z$ inducing the complex space topology. 
Assume that $g\colon P\times V \to Z_0$ is an 
$X$-holomorphic section of the form {\rm (\ref{eq:g})} 
and $G_t\colon P\times V \to Z$ $(t\in [0,1])$ is a 
homotopy of sections {\rm (\ref{eq:G})} satisfying 
\begin{itemize}
\item[\rm (a)] $\psi\circ G_t=g$ for $t\in [0,1]$,
\item[\rm (b)] $G_t(p,\cdotp)$ is holomorphic on $U$ for  
$(t,p) \in  Q$, and
\item[\rm (c)] 
$G_t(p,\cdotp)=G_0(p,\cdotp)$ for 
$(t,p) \in Q_0$, and these are holomorphic on $V$.
\end{itemize}
Let $\epsilon>0$. After shrinking the neighborhoods 
$U\supset K$ and $V\supset L$,  there exists a homotopy
$\wt G_t \colon P\times V \to Z$ $(t\in [0,1])$ 
of the form {\rm (\ref{eq:G})} such that 
\begin{itemize}
\item[\rm (i)]  $\psi\circ \wt G_t=g$ for all $t\in [0,1]$,
\item[\rm (ii)] for each $(t,p)\in Q$ the map 
$\wt G_t(p,\cdotp) \colon V \to Z$ is holomorphic 
and it satisfies 
$\sup_{x\in K} {\rm dist}\bigl(\wt G_t(p,x),G_t(p,x)\bigr) <\epsilon$,
and 
\item[\rm (iii)]  
$\wt G_t(p,\cdotp)=G_t(p,\cdotp)$ for each $(t,p) \in Q_0$.
\end{itemize}
Furthermore, there is a homotopy from $\{G_t\}$ to $\{\wt G_t\}$ 
consisting of homotopies with the same properties as $\{G_t\}$.
\end{proposition}

For families of sections of a holomorphic 
submersion $\pi\colon Z\to X$ over a Stein space $X$, PHAP holds if 
$Z\to X$ admits a fiber-dominating spray over a neighborhood 
of $L$ \cite{Gromov:OP,FPrezelj:OP2}, or a finite 
fiber-dominating family of sprays \cite{FF:subelliptic}.
Submersions admitting such sprays over small open subsets of $X$
are called {\em elliptic}, resp.\ {\em subelliptic}.
If PHAP holds over small open subsets of $X$ then sections $X\to Z$
satisfy the parametric Oka property (Gromov \cite[Theorem 4.5]{Gromov:OP};
the details can be found in \cite{FF:Kohn,FPrezelj:OP2}).
The same proof applies in our situation 
(see \cite[Theorem 4.2]{FF:Rothschild}) and shows 
that PHAP implies Theorem \ref{BOP-POP}. 
% Hence it remains to prove Proposition \ref{BOP-PHAP}.

%
%
%
%  THE PROOF
%
%
\section{Proof of Proposition \ref{BOP-PHAP}} 
Let $h\colon \cE\to Z$ denote the holomorphic vector bundle whose fiber
over a point $z\in Z$ equals the tangent space at $z$ to the (smooth)
fiber of $\psi$. 
The restriction $\cE|_\Omega$ to any open Stein domain $\Omega\subset Z$ 
is a reduced Stein space. By standard techniques we obtain for every 
such $\Omega$ an open Stein neighborhood $W \subset \cE|_\Omega$ of the zero 
section $\Omega \subset \cE|_\Omega$, with $W$ Runge in $\cE|_\Omega$,
and a continuous map
$s\colon \cE|_\Omega \to Z$ satisfying $\psi\circ s=\psi\circ h$,
such that $s$ is the identity on the zero section,
it is holomorphic on $W$, and it maps the fiber 
$W_{z}=\cE_z\cap W$ over a point $z\in Z$ biholomorphically 
onto a neighborhood of the point $z=s(0_z)$ 
in the fiber $Z_{\psi(z)}=\psi^{-1}(\psi(x))$.
Such $s$  is a fiber-dominating spray in the sense of 
\cite{Gromov:OP}, except that it is not globally holomorphic. 

\begin{figure}[ht]
\label{F1}
\psset{unit=0.5cm,linewidth=0.7pt}

\begin{pspicture}(-8.3,-0.5)(12,10)

\pspolygon(-3,0)(8,0)(12,6)(1,6)             %  Z
\psline[linewidth=0.4pt](-1,3)(10,3)         %  Z_x

\psline(-5,0)(-1,6)
\rput(-3.8,3){$z_0$}
\rput(-5.4,1.2){$Z_0$}

\psdot[dotsize=3pt](-3,3)
\psline[linewidth=0.2pt]{<-}(-2.8,3)(-1.2,3)
\rput(-1.8,3.5){$\psi$}

\pscurve[linewidth=1pt](0.5,0)(1,0.5)(3,2.5)(5,5)(5.8,6)       % f(X)
\psecurve[linewidth=1pt](0.5,0)(1,0.5)(3,2.5)(5,5)(5.8,6)      % f(V)
\psecurve[linewidth=1pt](2.5,0)(3,0.5)(5,2.5)(6.9,5)(7.5,6)    % g(V)

\pscustom[fillstyle=solid,fillcolor=lightgray]      % Z|_{f(U)}, shaded
{
\pscurve(0.5,0)(1,0.5)(3,2.5)(5,5)(5.8,6)
\psline(5.8,6)(5.8,10)
\pscurve(5.8,10)(5,9.3)(3,7)(1,5)(0.5,4.5)
\psline(0.5,4.5)(0.5,0)
}

\psecurve[linewidth=1pt](0.5,1.8)(1.1,2.3)(3,4)(5,6.7)(5.8,7.5)  % \xi_t 
\rput(2.7,5.2){$\xi_t$}
\psline[linewidth=0.2pt]{->}(2.7,4.7)(2.7,3.8)

\rput(9.6,4.7){$Z_{z_0}$}
\psline[linewidth=0.2pt]{->}(9.5,4.2)(9.5,3.05)

\rput(5.5,1.2){$G_t$}
\psline[linewidth=0.2pt]{<-}(3.8,1.2)(4.8,1.2)

\rput(9,7.8){$G_0(P\times V)$}
\psline[linewidth=0.2pt]{<-}(5.5,5.5)(8.2,7.3)

\rput(1,8.1){$\cE|_{G_0(P\times V)}$}
\psline[linewidth=0.2pt]{->}(1.3,7.5)(2.5,6.5)

\psarc[linewidth=0.2pt,arrows=<-](2,2){4}{21}{58}     %  curved arrow for $s$
\rput(5.6,4.5){$s_0$}

\rput(8,1.2){$Z$}

\end{pspicture}
\caption{Lifting sections $G_t$ to the spray bundle $\cE|_{G_0(P\times V)}$}
\label{Fig:1}
\end{figure}
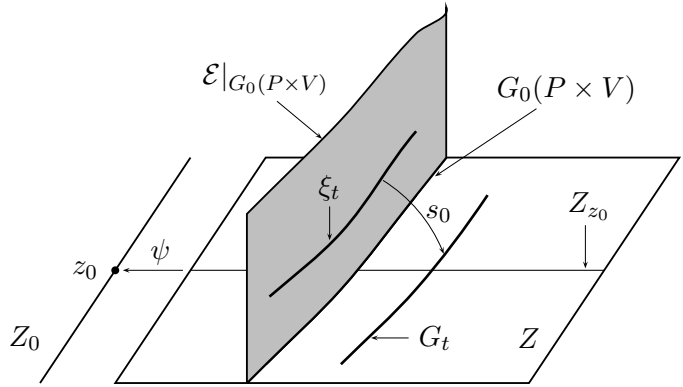

By \cite[Corollary 2.2]{F-Wold} each of the sets 
\[
		S_0= G_0(P\times L)\subset Z,\quad
		\Sigma_t = G_t(P\times K) \subset Z \quad (t\in [0,1])
\]
is a Stein compactum in $Z$.
By compactness of $\cup_{t\in[0,1]} \Sigma_t$
there exist numbers $0=t_0<t_1<\ldots <t_N=1$,
Stein domains $\Omega_0,\ldots, \Omega_{N-1}\subset Z$ 
satisfying 
\begin{equation}
\label{Omegaj}
		\bigcup_{t\in [t_{j},t_{j+1}]} \Sigma_t \ \subset \Omega_j
		\qquad (j=0,1,\ldots,N-1),
\end{equation}
and for every $j=0,1,\ldots,N-1$ there exist 
an open Stein neighborhood $W_j\subset \cE|_{\Omega_j}$ of the zero section
$\Omega_j$ of $\cE|_{\Omega_j}$ such that $W_j$ is Runge 
in $\cE|_{\Omega_j}$ and has convex fibers,
a fiber-spray $s_j \colon \cE|_{\Omega_j}\to Z$ as above 
that is holomorphic on $W_j$, 
and a homotopy $\xi_t$ $(t\in [t_{j},t_{j+1}])$ of 
$X$-holomorphic sections 
of the restricted bundle $\cE|_{G_{t_j}(P\times U)}$, 
with the range contained in $W_j$, such that 
\begin{itemize}
\item[(i)]  $\xi_{t_j}$ is the zero section of $\cE|_{G_{t_j}(P\times U)}$,  
\item[(ii)] $\xi_{t}(p,\cdotp)$ is the zero section when $p\in P_0$
and $t\in [t_j,t_{j+1}]$, and  
\item[(iii)] $s_j \circ \xi_t \circ G_{t_{j}} = G_t$ on $P\times U$
for all $t\in [t_j,t_{j+1}]$.
\end{itemize}
(See Fig.\ \ref{Fig:1}.)
For a given collection $(\Omega_j,W_j,s_j)$ the existence 
of homotopies $\xi_t$ is stable under sufficiently small 
perturbations of the homotopy $G_{t}$. 

%
%
%   LIFTING A SECTION TO A SPRAY BUNDLE
%
%

Consider the homotopy of sections $\{\xi_t\}_{t\in [0,t_1]}$ of 
$\cE|_{G_0(P\times U)}$. By the parametric version of 
the Oka-Weil theorem we can approximate 
$\xi_t$ uniformly on $P\times K$ by $X$-holomorphic sections 
$\wt \xi_t$ of $\cE|_{G_0(P\times V')}$ for an open set
$V'\subset X$ with $L\subset V'\subset V$. 
Further, we may choose $\wt\xi_t=\xi_t$ for $t=0$ and 
on $P_0\times V'$. In the sequel the set $V'$ may shrink around $L$.

By \cite[Corollary 2.2]{F-Wold} there is an open Stein neighborhood 
$\Omega$ of $S_0$ in $Z$ such that $S_0$ is $\cO(\Omega)$-convex. 
Hence $\Sigma_0=G_0(P\times K) \subset S_0$ 
is also $\cO(\Omega)$-convex, and it follows that
$W_0\cap \cE|_{\Sigma_0}$ is exhausted by 
$\cO(\cE|_\Omega)$-convex compact sets. 
Since $\cE|_{\Omega}$ is a reduced Stein space and $s_0$
extends continuously to $\cE|_\Omega$ preserving the
property $\psi\circ s_0=\psi\circ h$, the assumed 
BOP of $\psi$ implies that $s_0$ can be 
approximated on the range of the homotopy 
$\{\xi_t\colon t\in [0,t_1]\}$ 
(which is contained in $W_0\cap \cE|_{\Sigma_0}$) 
by a holomorphic map $\wt s_0\colon \cE|_{\Omega}\to Z$ 
which equals the identity on the zero section
and satisfies $\psi\circ \wt s_0=\psi\circ h$. 
The homotopy  
\[
		\wt G_t = \wt s_0 \circ \wt \xi_t\circ G_{0}
		\colon P\times V' \to Z  \quad (t\in [0,t_1])
\]
is fixed over $P_0$, $X$-holomorphic on $V'$, 
$\wt G_0=G_0$, and $\wt G_t$ approximates $G_t$ 
uniformly on $P\times K$ (also uniformly 
with respect to $t\in [0,t_1]$).
If the approximation is sufficiently close, 
we obtain a new homotopy $\{G_t\}_{t\in [0,1]}$ that agrees 
with $\wt G_t$ for $t\in [0,t_1]$ (hence is $X$-holomorphic on
$L$), and that agrees with the initial
homotopy for $t\in  [t'_1,1]$ for some $t'_1>t_1$ close to $t_1$.

We now repeat the same argument with the parameter interval 
$[t_1,t_2]$ using $G_{t_1}$ as the initial reference map. 
This produces a new homotopy that is $X$-holomorphic on $L$ 
for all values $t\in [0,t_2]$. 
After $N$ steps of this kind we obtain a homotopy 
satisfying the conclusion of Proposition \ref{BOP-PHAP}.

\medskip
\noindent\bf Acknowledgments. \rm  
I wish to thank Finnur L\'arusson
for his helpful remarks on a preliminary version of the paper.
Research was supported in part by grants P1-0291 and J1-2043-0101
from ARRS, Republic of Slovenia.
% , and by the Erwin Schr\"odinger Institut in Vienna, Austria.


\begin{thebibliography}{00}

\bibitem{FF:subelliptic}
Forstneri\v c, F.:
The Oka principle for sections of subelliptic submersions.
Math.\ Z.\ \textbf{241} (2002) 527--551 

\bibitem{FF:CR}
Forstneri\v c, F.: Oka Manifolds. 
C.\ R.\ Acad.\ Sci.\ Paris, Ser.\ I 347 (2009) 1017–-1020

\bibitem{FF:Kohn}
Forstneri\v c, F.:
The Oka principle for sections of stratified fiber bundles.
Pure and Appl.\ Math.\ Quarterly, 6 (2010), no.\ 3, 843--874


\bibitem{FF:Rothschild}
Forstneri\v c, F.:
Invariance of the parametric Oka property.
%Proceedings of the conference in honor of Linda P.\ Rothschild,
%Fribourg, 2008, 
P.\ Ebenfelt, N.\ Hungerbuehler, J.\ J.\ Kohn, 
N.\ Mok, E.\ J.\ Straube, eds., Complex  Analysis.  
Trends in Mathematics, Birkh\"auser (2010)


\bibitem{FPrezelj:OP2}
Forstneri\v c, F., Prezelj, J.:
Oka's principle for holomorphic submersions with sprays.
Math.\ Ann.\ 322 (2002) 633--666

\bibitem{F-Wold}
Forstneri\v c, F., Wold, E.\ F.:
Fibrations and Stein Neighborhoods.
Proc.\ Amer.\ Math.\ Soc., to appear. arXiv: 0906.2424 

\bibitem{Gromov:OP}
Gromov, M.:
Oka's principle for holomorphic sections of elliptic bundles.
J.\ Amer.\ Math.\ Soc.\ 2 (1989) 851-897

\bibitem{IK}
Ivarsson, B., Kutzschebauch, F.:
A solution of Gromov's Vaserstein problem.
C.\ R.\ Acad.\ Sci.\ Paris, Ser.\ I 346 (2008) 1239--1243

%\bibitem{IK2}
%Ivarsson, B., Kutzschebauch, F.:
%Holomorphic factorization of maps into $\mathrm{SL}_n(\C)$.
%Preprint (2008)

\bibitem{Larusson2}
L\'arusson, F.:
Model structures and the Oka principle.
J.\ Pure Appl.\ Algebra\ 192 (2004) 203--223 

\bibitem{Larusson3}
L\'arusson, F.:
Mapping cylinders and the Oka principle.
Indiana Univ.\ Math.\ J.\ 54  (2005) 1145--1159

\bibitem{Larusson4}
L\'arusson, F.: What is an Oka manifold?
Notices Amer.\ Math.\ Soc.\ 57 (2010), no.\ 1, 50--52.
http://www.ams.org/notices/201001/


\bibitem{Vaser}
Vaserstein, L.:
Reduction of a matrix depending on parameters to a diagonal 
form by addition operations.
Proc.\ Amer.\ Math.\ Soc.\ 103 (1988) 741--746 

\end{thebibliography}
\end{document}